\documentclass[11pt]{amsart}  
\usepackage{amsmath,amssymb,amsfonts,times,color,hyperref, fullpage}

\title[Classification of five dimensional DV-polyhedra]{The complete classification of five-dimensional Dirichlet-Voronoi polyhedra of translational lattices}

\author[M.~Dutour Sikiri\'c]{Mathieu Dutour Sikiri\'c}
\address{Rudjer Boskovi\'c Institute, Bijenicka 54, 10000 Zagreb, Croatia}
\email{mdsikir@irb.hr}

\author[A.~Garber]{Alexey Garber} 
\address{School of Mathematical \& Statistical Sciences, The University of Texas Rio Grande Valley, 1 West University Blvd, Brownsville, TX, 78520, USA.}
\email{alexeygarber@gmail.com}

\author[A.~Sch\"urmann]{Achill Sch\"urmann} 
\address{Institute of Mathematics, University of Rostock, 18051 Rostock, Germany}
\email{achill.schuermann@uni-rostock.de}

\author[C.~Waldmann]{Clara Waldmann}
\address{Department of Mathematics, Technical University Munich, Boltzmannstr.3, 85747 Garching, Germany}
\email{clara.waldmann@tum.de}


\date{\today}
\newtheorem{defi}{Definition}[section]

\newtheorem{theorem}[defi]{Theorem}

\newcommand{\Z}{{\mathbb{Z}}} 
\newcommand{\R}{{\mathbb{R}}} 
 
\newcommand{\MD}{{\mathcal{D}}} 
\newcommand{\MS}{{\mathcal{S}}}

\DeclareMathOperator{\DV}{DV}
\DeclareMathOperator{\conv}{conv}
\DeclareMathOperator{\GL}{GL}
\DeclareMathOperator{\Del}{Del}
\DeclareMathOperator{\SC}{SC}
\DeclareMathOperator{\Can}{Can}
\DeclareMathOperator{\Trace}{Trace}

\renewcommand{\proof}{\noindent{\bf Proof.}\ \ }
\renewcommand{\qed}{\hfill $\Box$ }

\newcommand{\sd}{{\mathcal S}^{d}} 
 
\newcommand{\sdgo}{{\mathcal S}^{d}_{>0}} 
\newcommand{\sdclosed}{{\mathcal S}^{d}_{rat,\geq 0}}
\newcommand{\Stab}{\operatorname{Stab}}
\newcommand{\id}{\operatorname{Id}}

\begin{document}

\begin{abstract} 
In this paper we report on the full classification of Dirichlet-Voronoi polyhedra 
and Delaunay subdivisions of five-dimensional translational lattices. 
We obtain a complete list of $110244$ affine types (L-types)
of Delaunay subdivisions 
and it turns out that they are all combinatorially inequivalent,
giving the same number of combinatorial types of Dirichlet-Voronoi polyhedra.
Using a refinement of corresponding secondary cones, we obtain 
$181394$ contraction types.
We report on details of our computer assisted enumeration, which we
verified by three independent implementations and a topological mass formula check.
\end{abstract}

\maketitle


\section{Introduction}

\sloppy

The study of translational lattices and their Dirichlet-Voronoi polyhedra are classical subjects
in crystallography.
In 1885 Fedorov \cite{fedorov-1885} (cf. \cite{sg-1984})
determined the five combinatorial types of possible Dirichlet-Voronoi polyhedra in 
the Euclidean $3$-space $\R^3$.
These are also all the parallelohedra in $\R^3$, that is, 
polyhedra admitting a facet-to-facet tiling of $\R^3$ by translation.
Voronoi~\cite{voronoi-1908} developed a theory 
to classify Dirichlet-Voronoi polyhedra
for arbitrary $d$-dimensional Euclidean spaces~$\R^d$.
His theory allows to classify them via a classification of Delaunay subdivisions 
up to affine equivalence (so called {\em L-types}).
In this context Voronoi also 
came up with his famous and still unsolved conjecture, stating that
every parallelohedron in $\R^d$ is affinely equivalent to a 
Dirichlet-Voronoi polyhedron for some translational lattice.

In this paper we report on the enumeration of the $5$-dimensional combinatorial types of
Dirichlet-Voronoi polyhedra or equivalently Delaunay subdivisions
(Theorem~\ref{thm:ourmainthm-forDelone}). 
We find in total $110244$ different combinatorial types
and hereby go beyond the partial classification according to
subordination schemes previously obtained by \cite{engel-2000}.
In Table~\ref{tab:NumberLtypeEnumeration} we list the number
of Delaunay subdivisions that were computed so far. 
By our work, a full classification is known for $d\leq 5$ so far.
Recent partial results on primitive types in dimension~$6$
\cite{be-2013} seem to indicate that a full classification beyond 
$5$~dimensions is out of reach at the moment.

Our paper is organized as follows.
In Section~\ref{sec:d-polytopes} we start with some notation and background 
on Dirichlet-Voronoi and Delaunay polytopes.
Voronoi's L-type theory
is briefly reviewed in Section~\ref{sec:voro-theory}.
We in particular describe how the classification of Dirichlet-Voronoi polyhedra 
is reduced to the classification of Delaunay subdivisions and how 
this can practically be done.
Algorithms and implementations for our classification result 
are briefly described in Section~\ref{sec:algorithms_and_implementations}
and references to online sources are given.
Additional data and tables are presented in
Section~\ref{sec:tables_and_data},
where we also relate our work to the theory of contraction types.

\section{Dirichlet-Voronoi and Delaunay polytopes}  \label{sec:d-polytopes}

Let $\Lambda$ denote a {\em translational lattice} in $\R^d$.
That is, $\Lambda$ is a full rank-discrete subgroup of $\R^d$
and, equivalently, can be written as
\begin{equation*}
\Lambda =
\left\{
\lambda_1 b_1 +\ldots +\lambda_d b_d
\; : \;
\lambda_1,\ldots,\lambda_d \in \Z
\right\}
\end{equation*}
with linearly independent vectors $b_1,\ldots, b_d \in \R^d$.
Latter vectors, as well as a matrix $B$ with these as columns, 
are referred to as a {\em basis} of $\Lambda$ and
we simply write $\Lambda=B\Z^d$.
Viewing $\R^d$ as a Euclidean space with norm
$|\cdot|$, the {\em Dirichlet-Voronoi polytope} of $\Lambda$ is
defined as the set of points in $\R^d$ which are at least as close to the origin
than to any other element of $\Lambda$: 
\begin{equation*}
\DV(\Lambda) = \left\{
x\in \R^d\; : \;
|x| \leq |x-y| \quad \mbox{for all } y\in \Lambda
\right\}.
\end{equation*}

\subsection{General facts about polytopes} \label{sec:polytope-background}

The term polytope refers to the fact that $\DV(\Lambda)$ can
be described as a convex hull (set of all convex combinations) of
finitely many points. A point that can not be omitted in such a description
is called a {\em vertex} of the polytope.
Let us briefly review some basics from the theory of polytopes
(see \cite{ziegler-1997,gruenbaum-2003} for details).
A {\em supporting hyperplane} is an affine
hyperplane having the property that the polytope is fully
contained in one of the two halfspaces bounded by it.
A {\em $k$-dimensional face} of a polytope is defined as
a $k$-dimensional intersection of the polytope with a
supporting hyperplane.
The $(d-1)$-dimensional faces of a $d$-dimensional polytope are called {\em facets}
and vertices are the $0$-dimensional faces.
Every polytope also has a description by linear inequalities and
the non-redundant ones in such a description are 
in 1-to-1-correspondence to its facets.

Altogether, the faces of a polytope 
form a poset (partially ordered set, ordered by inclusion),
which is called the {\em face lattice} of the polytope.
Two polytopes are called {\em combinatorially equivalent}, 
if they possess the same face lattice.
For instance, two $2$-dimensional $n$-gons 
(which are the $2$-dimensional polytopes with $n$ vertices)
are always combinatorially equivalent.
However, they might not be {\em affinely equivalent}, that is, there does not exist an
affine map, mapping one to the other
(see \cite{GroupPolytopeLMS} for details on this
and how to compute equivalence).

We note that Engel \cite{engel-2000} uses 
a so called {\em subordination scheme} 
(sometimes called {\em polyhedral scheme})
which is an invariant to classify Dirichlet-Voronoi polytopes.
Two combinatorially different polytopes can however have the
same subordination scheme. In fact, several 
combinatorially different Dirichlet-Voronoi polyhedra in $\R^5$
have the same subordination scheme.
Therefore this invariant can not be used for a full
classification of all combinatorial types.

\subsection{Affine and combinatorial types of Dirichlet-Voronoi polytopes}

In dimension~$2$ there exist only two combinatorially inequivalent types
of Dirichlet-Voronoi polytopes: either centrally symmetric hexagons or rectangles.
We note that there are infinitely many affine types of Dirichlet-Voronoi polytopes.
Actually, any centrally symmetric hexagon with vertices on a unit circle
is a Dirichlet-Voronoi polytope of a lattice.
However, they are not all affinely equivalent to each other.
For instance, none of them is affinely equivalent to a regular hexagon 
(except the regular hexagon itself).
We refer the interested reader for more information on affine types
of Dirichlet-Voronoi polytopes to~\cite{din-2011,gavrilyuk-2014}.

The combinatorial types of Dirichlet-Voronoi polytopes in dimensions~$3$ and~$4$
are known as well.
There exist five different combinatorial types of Dirichlet-Voronoi polytopes
in dimension~$3$ and $52$ different combinatorial types in dimension~$4$.
In this paper we report on the classification in dimension~$5$ 
and we show:

\begin{theorem}\label{thm:ourmainthm}
There are precisely 
$110244$ combinatorially inequivalent types of Dirichlet-Voronoi polytopes
of five-dimensional translational lattices.
\end{theorem}

In the following we explain in more detail how 
to obtain the above classification result, based on Voronoi's second reduction theory
for positive definite quadratic forms.

\subsection{Delaunay subdivisions}

The notion of Delaunay subdivisions were introduced in \cite{delone-1934}.
Here we give their definition and shortly describe major properties.

Given a translational lattice $\Lambda$ in $\R^d$,
an empty sphere $S(c,r)$ of center $c$ and radius $r>0$ is a sphere
such that there is no lattice point in its interior.
A {\em Delaunay cell} is an intersection $\Lambda \cap S(c,r)$.
A {\em Delaunay polytope} is a $d$-dimensional polytope of the
form $\conv(\Lambda\cap S(c,r))$.

The set of all Delaunay polytopes of $\Lambda$ form a
polytopal subdivision of~${\R}^d$, called the {\em Delaunay subdivision} of $\Lambda$.
In general, a polytopal subdivision is a non-overlapping union of polytopes
that fill all of~$\R^d$ and such that the intersection of any two polytopes is either empty or
a $k$-dimensional face.
$\DV(\Lambda)$ together with all its translates by lattice vectors
form another polytopal subdivision of~$\R^d$.
Both subdivisions are invariant by lattice translations.
The Delaunay polytopes with vertex at $x\in \Lambda$ are translates 
by $x$ of some Delaunay polytope with vertex at $0$. Thus to know the 
full Delaunay subdivision of a lattice $\Lambda$, it suffices to know 
the Delaunay polytopes with vertex~$0$. 
The centers of these Delaunay polytopes coincide with the vertices of $\DV(\Lambda)$.

The Delaunay subdivision is said to be {\em dual} to the subdivision with 
Dirichlet-Voronoi polytopes. The Dirichlet-Voronoi polytope of a lattice can be
obtained from the Delaunay polytopes with vertex~$0$ and vice versa:
There is a bijection between the $k$-dimensional faces of these Delaunay polytopes
and the $(d-k)$-dimensional faces of the Dirichlet-Voronoi polytope.
In particular, each $d$-dimensional Delaunay polytope corresponds to a vertex 
of the Dirichlet-Voronoi polytope.
Moreover, the face lattice structure with respect to inclusion is preserved as well:
If two faces of Delaunay polytopes with vertex~$0$ 
are contained in each other,
the corresponding dual faces of the Dirichlet-Voronoi polytope 
are contained in each other with the inclusion reversed.
Therefore, the classification of combinatorial types of Dirichlet-Voronoi polytopes is 
equivalent to the classification of combinatorial types of Delaunay subdivisons.

The different combinatorial types can be derived from possible affine types.
Here, two Delaunay subdivisons, respectively lattices $\Lambda$ and $\Lambda'$ 
are {\em affinely equivalent} (are of the {\em same affine type}),
if there is a matrix (linear map) $A\in\GL_d(\R)$ with $\Lambda'=A\Lambda$, 
mapping all Delaunay polytopes of $\Lambda$ to those of $\Lambda'$.
Note that two Delaunay subdivisions 
with different combinatorial types can not be affinely equivalent.
The opposite could be possible though:
Two different affine types of Delaunay subdivisions could 
possibly have the same combinatorial type --- although we do
not know of a single example among Delaunay subdivisons 
for translational lattices at this point.
In particular, up to dimension~$5$, all 
affine types of Delaunay subdivisons are not only affinely
inequivalent, but also combinatorially.

\section{Voronoi's second reduction theory}  \label{sec:voro-theory}

In the following we give a short sketch of Voronoi's second reduction theory
\cite{voronoi-1908}, as far as it is necessary to describe how our classification of
affine types of five-dimensional Delaunay subdivisions is obtained.
For a more detailed description and extensions 
of the theory we refer to~\cite{schurmann-2009}.

\subsection{Working with Gram matrices}

The set of real symmetric positive definite matrices is denoted $\sdgo$.
When dealing with lattices up to orthogonal transformations, it is often
convenient to work with Gram matrices $Q=B^tB\in\sdgo$
instead of using matrices of lattice bases~$B$. 
Up to orthogonal transformations, the basis matrix $B$ can uniquely
be recovered from $Q$ using the Cholesky decomposition.
Geometrically this is equivalent to reconstruction 
of a basis knowing vector lengths and angles between them.
Every positive definite symmetric matrix~$Q$ defines 
a corresponding positive definite quadratic form $x\mapsto Q[x] = x^t Q x$ on $\R^d$.

In particular for studying affine types of Delaunay subdivisions it is convenient to
use the same coordinates of vertices $v_1,\ldots, v_n$ from a fixed 
translational lattice $\Lambda\subseteq \R^d$ (often $\Lambda= \Z^d$)
for different affine images $B\cdot \conv \{v_1,\ldots, v_n\}$ of Delaunay polytopes, 
which we represent by a corresponding matrix $Q\in\sdgo$.
A polytope $P = \conv\{v_1, \ldots , v_n \}$ with vertices $v_i \in \Lambda$ 
is called a {\em Delaunay polytope of~$Q$} if it is $d$-dimensional and 
if there exists a center $c\in \R^d$ and a real number $r$ such that
$Q [c-v_i] = r^2$ for $i=1\ldots , n$ and 
$Q [c- v] > r^2$ for all other $v\in\Lambda$. 
The set $\Del(\Lambda,Q)$ of all Delaunay polytopes of $Q\in\sdgo$
is a polytopal subdivision of $\R^d$,
called the {\em Delaunay subdivision of~$Q$ with respect to $\Lambda$}.

We speak of a {\em  Delaunay triangulation}, if all the Delaunay polytopes are
simplices, that is, if all of them have affinely independent vertices.
We say that $\Del(\Lambda , Q)$ is a {\em refinement}  of $\Del(\Lambda , Q')$ 
(and $\Del(\Lambda , Q')$ is a {\em coarsening} of $\Del(\Lambda , Q)$), 
if every Delaunay polytope of~$Q$ is contained in a Delaunay polytope of~$Q'$.
Any Delaunay subdivision can be refined to a Delaunay triangulation
by perturbing $Q$ if necessary.
Voronoi's theory of secondary cones which we explain below 
gives us an explicit description of the 
set of positive definite matrices having the same Delaunay subdivision.

\subsection{Secondary Cones and L-types}

Voronoi's second reduction theory is based on {\em secondary cones}
(also called {\em $L$-type domains})
\[
\SC(\MD) = \left\{Q \in \sdgo : \Del(\Z^d, Q) = \MD\right\},
\]
which can be seen to be non-empty polyhedral cones in $\sdgo$
(which are open within their linear hull),
if $\MD$ is a Delaunay subdivision for some~$Q$.
In order to give an explicit description of $\SC(\MD)$ we 
define for an affinely independent set $V\subseteq  \Z^d$ 
of cardinality $d+1$ and a point
$w\in\Z^d$ the symmetric matrix
\begin{equation} \label{eqn:N-forms}
N_{V, w} = w w^t - \sum_{v \in V}
\alpha_{v} v v^t
,
\end{equation}
where the coefficients $\alpha_{v}$ are uniquely determined by
the affine dependency 
$$
w = \sum_{v \in V} \alpha_{v} v 
\quad \mbox{with} \quad
1 = \sum_{v \in V} \alpha_{v}
.
$$
In the special situation of 
$V=\{v_1, \ldots,v_{d+1}\}$ being vertices of a Delaunay simplex~$L$
and $w$ being the additional vertex of a Delaunay simplex 
$L' = \conv\{v_{2}, \ldots, v_{d+1},w\}$ 
{\em adjacent}  to $L$, we use the notation $N_{L,L'}$ for $N_{V,w}$. 
In the following we use $\langle A, B \rangle = \Trace \left(AB\right)$ to denote the 
standard inner product defined for two symmetric matrices $A,B$ on $\sd$.
The following result by Voronoi gives an explicit description of
a secondary cone in terms of linear inequalities.

\begin{theorem}[\cite{voronoi-1908}]  
\label{thm:explicit-description-secondary-cones}
Let $Q$ be a positive definite symmetric matrix
whose Delaunay subdivision $\MD=\Del(\Z^d, Q)$ is a triangulation. 
Then 
\begin{equation}   \label{eqn:explicit-description-secondary-cones}
\SC(\MD) = \{Q' \in \sd :
\mbox{$\langle N_{L,L'}, Q'\rangle >0$
for adj. $L,L'\in \MD$}\}.
\end{equation}
\end{theorem}

This theorem of Voronoi shows that the secondary cone $\SC(\MD)$
of a Delaunay triangulation $\MD$ is a full dimensional {\em open polyhedral cone},
that is, the intersection of finitely many open halfspaces.
If we use weak inequalities $\geq 0$  in~\eqref{eqn:explicit-description-secondary-cones}
instead of strict inequalities, we obtain a description of the closed polyhedral cone $\overline{\SC(\MD)}$.
We will use these closed versions and their facial structure in the sequel.
Just like for polytopes (cf. Section~\ref{sec:polytope-background}), faces can be defined 
for these closed polyhedral cones and the set of all faces forms a combinatorial lattice -- 
{\em the face-lattice of the cone}.
Voronoi discovered that the faces of $\overline{\SC(\MD)}$ 
correspond to all the possible coarsenings of $\MD$.

Two full dimensional secondary cones touch in a facet, if and only if the corresponding
Delaunay triangulations can be transformed into each other by
{\em bistellar flips}.
That is we first apply a coarsening of some of the simplices to repartitioning
polytopes ($d$-dimensional polytopes with $d+2$ vertices) and then apply a
refinement procedure.
Since these changes of Delaunay triangulations are 
not important for what follows, we omit a detailed description here
and refer the interested reader to~\cite{schurmann-2009}.

The rational closure $\sdclosed$ of $\sdgo$ is the set of positive
semidefinite quadratic forms whose kernel is defined by rational equations.
At the core of Voronoi's theory is the action of the matrix group $\GL_d(\Z)$ on
the polyhedral tiling by closures of secondary cones:

\begin{theorem}[Voronoi's Second Reduction Theory]
\label{thm:mainvoronoi}
The topological closures $\overline{\SC(\MD)}$ give a polyhedral
subdivision of $\sdclosed$ on which the group $\GL_d(\Z)$ acts
by $\overline{\SC(\MD)} \mapsto U^t \overline{\SC(\MD)}U$. 
Under this group action there are only finitely many inequivalent
secondary cones.
\end{theorem}

Note that one can subdivide 
the secondary cones into smaller cones and obtain a reduction domain
for the action of ~$\GL_d(\Z)$ on $\sdgo$.
This is the reason why Voronoi's theory of Delaunay subdivisions and secondary cones
is referred to as {\em Voronoi's second reduction theory} (for positive definite quadratic forms).

For our classification of affine types 
the following observation is crucial:

\begin{theorem} \label{thm:affinetypes} 
Let $Q, Q'\in \sdgo$ be two positive definite matrices
with Cholesky decompositions $Q=B^t B$ and $Q'=(B')^t(B')$
and corresponding lattices $\Lambda=B\Z^d$ and $\Lambda'=B'\Z^d$.
Then the Delaunay subdivisons of $\Lambda$ and $\Lambda'$ are of the same
affine type if and only if $Q$ and $Q'$ are in $\GL_d(\Z)$-equivalent secondary cones.
\end{theorem}
\proof We are not aware of an explicit reference for this result, so for clarity
we give an argument here. First we note that transforming a set $\Lambda$
and a Delaunay decomposition $\Del(\Lambda, Q)$ 
by a linear map $A\in\GL_d(\R)$ we get a new 
Delaunay decomposition 
$\Del(\Lambda', (A^{-1})^tQA^{-1})$ with vertex set $\Lambda'=A\Lambda$.

Suppose now that the Delaunay decompositions of 
$\Lambda$ and $\Lambda'$ are of the same affine type.
Then
$A\cdot \Del(\Lambda, \id_d) = \Del(\Lambda', (A^{-1})^tA^{-1})  = \Del(\Lambda',\id_d)$.
Therefore 
\begin{equation*}
\begin{array}{rcl}
\Del(\Z^d, Q)
&=& B^{-1} \Del(\Lambda, \id_d)\\
&=& B^{-1} A^{-1} \Del(\Lambda', \id_d)\\
&=& U \Del(\Z^d, Q')
\end{array}
\end{equation*}
with 
$U= B^{-1} A^{-1} B'$.
Since $\Z^d = U \Z^d$ we have $U\in \GL_d(\Z)$
and therefore $Q$ and $(U^{-1})^t Q' U^{-1}$ are in the same secondary cone.

On the other hand, if $Q$ and $Q'$ are in 
$\GL_d(\Z)$-equivalent secondary cones, then 
there exists a $U\in \GL_d(\Z)$
with $\Del(\Z^d, Q') = U \Del(\Z^d, Q)$.
Thus
$$(B')^{-1}\Del(\Lambda', \id_d)= U B^{-1}\Del(\Lambda, \id_d),$$
and hence $A= B' U B^{-1}$ satisfies
$A \Del(\Lambda, \id_d) = \Del(\Lambda', \id_d)$. \qed

With the knowledge on how to perform bistellar flips,
Theorems~\ref{thm:mainvoronoi} and~\ref{thm:affinetypes} easily lead to an algorithm 
to enumerate all affine types of Delaunay triangulations in a given dimension
(see Algorithm 3 in \cite{schurmann-2009}). 
For it, Sch\"urmann and Vallentin developed 
the program {\tt scc} ({\em secondary cone cruiser}).
Its first version from \cite{scc} already allowed to reproduce
the known classification of all $\GL_d(\Z)$-inequivalent
Delaunay triangulations up to dimension~$d=5$. 
We will use their result, respectively the output of the program {\tt scc}. 

Beginning with dimension~$6$
the number of inequivalent Delaunay triangulations starts to explode.   
At the moment, we still do not know how  
many inequivalent triangulations we have to expect in dimension~$6$.
Baburin and Engel \cite{be-2013} report that they found 
$567.613.632$ ones so far.

\subsection{Enumeration of all Delaunay subdivisions} 
\label{sec:classification-algorithms}

Arbitrary Delaunay subdivisions are limiting cases of Delaunay triangulations. 
Their secondary cones occur on the boundaries of full-dimensional secondary
cones of Delaunay triangulations. 
The following theorem seems to be folklore.
One can find a proof for example in Proposition~2.6.1 of \cite{vallentin-2003}:

\begin{theorem}
\label{thm:delonedegeneration}
Let~$\MD$ be a Delaunay triangulation.
\begin{enumerate}
\item A positive definite symmetric matrix $Q$ lies
in~$\overline{\SC(\MD)}$ if and only if~$\MD$ is a refinement of~$\Del(Q)$.
\item If two positive definite symmetric matrices $Q$ and $Q'$ both
lie in~$\overline{\SC(\MD)}$, then $\Del(Q+Q')$ is a common refinement
of~$\Del(Q)$ and $\Del(Q')$.
\end{enumerate}
\end{theorem}

We note that this theorem can be extended to positive semidefinite symmetric matrices
in the rational closure $\sdclosed$ of $\sdgo$. For those among them which are 
not positive definite, one can define a polyhedral Delaunay subdivision 
with unbounded polyhedra.
For details we refer to Chapter~4 of \cite{schurmann-2009}.

By Theorem~\ref{thm:delonedegeneration},
the classification of all inequivalent Delaunay subdivisions
is equivalent to the classification of all inequivalent secondary cones.
In order to prove our Theorem~\ref{thm:ourmainthm},
we show the following equivalent result: 

\begin{theorem}  \label{thm:ourmainthm-forDelone}
In dimension~$5$ there are  
$110244$ affine types of Delaunay subdivisons.
Equivalently, there are that many secondary cones 
of positive definite quadratic matrices in $\MS^5$  
up to $\GL_5(\Z)$-equivalence.
\end{theorem}

\subsection{Related works} 
\label{sec:related-work}

At this point, we should point out that there is a parallel theory that considers
a single Delaunay polytope in a lattice, irrespective of the other Delaunay polytopes
in the tessellation. This theory is exposed in \cite{DL} and recent developments can be
found in \cite{InhomogeneousPerfect}.
The possible Delaunay polytopes of dimension~$5$ were classified in~\cite{kononenko-2002}
in terms of $138$ combinatorial types. The classification in dimension $6$ in \cite{DelaunaySix}
gives $6241$ combinatorial types.

In \cite{schurmann-2009} (cf. Table~2 on page~60) it is reported that 
Engel~\cite{engel-2000} found $179372$ inequivalent
five-dimensional Delaunay subdivisions.
This, however, is unfortunately a misinterpretation of Engel's result
who classifies so called {\em contraction types} (of parallelohedra).
From these contraction types, he derives $103769$ ``combinatorial types''.
These types are not the true combinatorial types that are classified here however,
but a coarser notion, which classifies parallelohedra in dimension~$5$, 
or equivalently Delaunay subdivisions, up to their subordination schemes.
The subordination scheme of a $d$-dimensional
polytope~$P$ is a list of numbers, containing 
for every $k=2,\ldots,d-1$ and for every $n$, 
the number of $(k-1)$-faces of~$P$ incident to exactly~$n$ of the 
$k$-faces of~$P$ (see Section~4 of \cite{engel-2000} for details). 
Thus, the subordination scheme encodes certain properties of the face lattice of
a polytope, but not the whole face lattice.
Two combinatorially different polytopes can have the same
subordination scheme. They may even be the same for different affine types of
Dirichlet-Voronoi polytopes, having even secondary cones of different
dimension. In fact, during our work we discovered two such examples
for $d=5$.

Note that combinatorial types of polytopes can only truely be
distinguished, by checking whether or not their face lattices are
different.
It has been shown in~\cite{KaibelSchwartz2003} that the incidence 
relations between vertices and facets of two polytopes
are sufficient to distinguish their face lattices. 
Practically such differences can be checked 
using graph isomorphism software as we describe in the next section.
Invariants like the number of faces of a given dimension or
the subordination scheme used by Engel may be useful in computations, 
for instance when limiting the number
of equivalence tests.
However, such invariants are not sufficient for complete enumerations.
Engel's invariant appears to distinguish the known~$52$ 
combinatorial types in dimension~$4$, but it does not distinguish 
types in any dimension greater or equal to~$5$.
While it is conceivable that the subordination
scheme could be extended to better distinguish between types, 
it should never be used alone without checking for equivalence since
there is always the possibility that non-isomorphic structures
have the same invariant.

\section{Algorithms and Implementations}  \label{sec:algorithms_and_implementations}

Before we explain the details of our computations for $d=5$, we
start with some general observations, which are valid in all dimensions
and quite useful for practical purposes.

\subsection{Using reduced generators and central forms}

Each secondary cone, respectively its closure is given by a finite list of linear inequalities 
(coming from Voronoi's regulators, cf. Theorem~\ref{thm:explicit-description-secondary-cones}). From it one can obtain a number of {\em generating rays}.
In fact, one of these descriptions (by rays or inequalities) can be obtained from the other by
a polyhedral representation conversion.
Since all of the involved inequalities involve rational numbers only, we 
may assume that the generators for rays are given by integral vectors (matrices in $\sd$), 
with coordinates having a greatest common divisor (gcd) of~$1$.
We refer to these generators as {\em reduced (or normalized) generators}.
As we are using Theorem~\ref{thm:delonedegeneration} for the
classification of Delaunay subdivisons, we only need to consider
closures of secondary cones which are faces of
closures of full-dimensional secondary cones. 
All such faces are themselves generated by a subset of the 
reduced generators of the full dimensional cone.

Having reduced generators $R_1,\ldots, R_k$ of a secondary cone $\SC$ 
(respectively its closure), 
we define a {\em central reduced (or normalized) form} of the secondary cone as the sum
$Q(\SC)=\sum_{i=1}^k R_i$.
It is easy to see that two secondary cones $\SC$ and $\SC'$ are $\GL_d(\Z)$-equivalent
if and only if $Q(\SC)$ and $Q(\SC')$ are $\GL_d(\Z)$-equivalent.
Hence, for the classification of secondary cones up to $\GL_d(\Z)$-equivalence 
we can equally well classify their central reduced forms up to $\GL_d(\Z)$-equivalence.

\subsection{Testing equivalence of forms and use of invariants}

Testing $\GL_d(\Z)$-equivalence 
of central reduced forms can be done with the Plesken-Souvignier algorithm \cite{plesken1993}.
Their initial implementation is available at \cite{isom} and is part of 
computer algebra software such as \cite{magma,gap}.
The algorithm works by building a finite set of vectors that is
canonically defined by a given positive definite matrix and spans $\Z^d$ as a lattice.
For a given norm bound $n$ and a positive definite matrix $Q$ let
\begin{equation*}
  S(Q,n) = \left\{ v\in \Z^d \mbox{~s.t.~} Q[v] \leq n\right\}.
\end{equation*}
Then we take the smallest $n$ such that $S(Q,n)$ spans $\Z^d$ as a lattice and call the
vector set $\Can(Q)$.


As testing $\GL_d(\Z)$-equivalence of central reduced forms is computationally quite 
involved, one needs to reduce the number of such tests as much as possible since the final
number of forms is $M=110244$ and so the total number of isomorphism tests is a priori
$M(M-1)/2$.
The basic idea is to use invariants to reduce the number of tests.
Some invariants come naturally from the form $Q(\SC)$ such as its determinant and
size of $\Can(Q(\SC))$.
Other possible invariants are related to the secondary cone $\SC$ under consideration.
For example the dimension of $\SC$ or its number of generating forms $R_1$, \dots, $R_k$.
Further invariants are the rank of $R_k$ and so on.
Rather surprisingly, the most efficient invariant tends to be the determinant of $Q(\SC)$.

\subsection{Putting it all together for five dimensions}

Now, finally, let us put the pieces above together, to describe the
algorithm behind our classification result for $d=5$.
To show Theorem \ref{thm:ourmainthm-forDelone} 
with computer assistance, we can use Voronoi's theory.
We start from the secondary cones of the $222$
known Delaunay triangulations.
Those were classified in \cite{br-1973,rb-1976}
but the classification was incorrect and a final correct classification
was obtained in \cite{eg-2002} which we have independently confirmed
in \cite{CompApproachesLatticeCovering,DecompHypermetricCone}.
These open polyhedral cones are full dimensional in ${\mathcal S}^5_{>0}$ 
and therefore have dimension~$15$.
Their closure is given by a list of non-redundant linear inequalities.
From this list, we can obtain the reduced generators of each cone
and also a description by generators 
and by equations / inequalities for each of their facets.
These facets are themselves closures of 
$14$-dimensional secondary cones which 
correspond to Delaunay subdivisions that are a true coarsening of the considered  
Delaunay triangulation at hand. 
Some of them may be $\GL_d(\Z)$-equivalent, so for our classification,
we have to obtain a list of  $\GL_d(\Z)$-inequivalent $14$-dimensional secondary cones 
in ${\mathcal S}^5_{>0}$ 
from them, using their central reduced forms.
In a next step, we obtain a list of  $\GL_d(\Z)$-inequivalent $13$-dimensional 
secondary cones from our list of $14$-dimensional secondary cones in a similar way.
We continue this process until we subsequently obtain a full list of 
$\GL_d(\Z)$-inequivalent cones of dimensions $15, \ldots, 1$.
See Table~\ref{tab:number-secondary-cones-contraction-cone-bydim} for the number of secondary cones 
obtained in each dimension in this way.

\subsection{Practical Implementations}\label{sec:classification-results}

The computer code of our first implementation in {\tt Haskell}
of the algorithm described above, together with a detailed documentation (in German)
is available at the webpage~\cite{waldmann}.
In particular, data of the full classification can be obtained at~\cite{fulldata5}, 
with a matrix of a central reduced form 
for each secondary cone in ${\mathcal S}^5_{>0}$.

Our second implementation used the {\tt GAP} package {\tt polyhedral}
\cite{polyhedral} 
with some external calls to {\tt isom} \cite{isom} for equivalence tests
and {\tt lrs} \cite{lrs} for polyhedral representation conversions.
In our third implementation, we adapted the program {\tt scc}.
In its latest version \cite{scc2} we  
included the program {\tt isom}  
to produce all secondary cones of a given dimension.

In order to avoid the dependency on {\tt isom} 
in all three implementations, we also performed equivalence computations with
{\tt nauty} \cite{nauty}, applied to test equivalence of the sets $\Can(Q(\SC))$ of vectors,
by using the method explained in Section~3.4 of \cite{GroupPolytopeLMS}.
Overall, the full computation, its resulting data and in particular the numbers in
Table~\ref{tab:number-secondary-cones-contraction-cone-bydim} were all sufficiently
well cross-checked.
All calculations yield the same results
and due to the different nature of our three programs
we can be certain of the obtained classification,
although the computations are large and quite involved.

We can use the obtained results for a computational proof of our main 
Theorem~\ref{thm:ourmainthm}, 
by showing that all
Delaunay subdivisions, respectively the corresponding Dirichlet-Voronoi polytopes,
are combinatorially inequivalent.
This is shown by checking if their face-lattices are non-isomorphic.
Since the face-lattice of a polytope is determined by the incidence graph
of vertices and facets,
we can check if these graphs are non-isomorphic.
These isomorphism checks can be performed using for instance 
graph isomorphism software like~{\tt nauty} \cite{nauty}.
We computed ``canonical forms'' for each of the graphs with 
{\tt nauty} 
and then used {\tt md5sum} (a special hash function) for each of them
in order to decide computationally (in a reasonable amount of time)
that they are all different.

\section{Tables and Data}\label{sec:tables_and_data}

We provide the following tables, containing additional information:
Table \ref{tab:number-secondary-cones-contraction-cone-bydim}
gives the number of inequivalent secondary cones by their dimension.
Table \ref{tab:number-secondary-cones-byrank} gives the number of secondary
cones by their number of rank $1$, $4$ or $5$ extreme rays.
Table \ref{tab:NumberLtypeEnumeration} gives the known numbers
of inequivalent secondary cones (all combinatorial types) 
and full-dimensional secondary cones (primitive types), 
together with a reference where these results can be found.
Table \ref{tab:number-domain-by-dim-and-gen} gives the number of secondary
cones according to their dimension and their number of extreme rays.
Table \ref{tab:number-secondary-cones-not-extendable} gives the number of
secondary cones that cannot be extended to a higher dimensional cone by
a pyramid construction with a rank-$1$ extreme ray.
Table \ref{tab:freq-bravais-group} gives the frequencies of occurring
Bravais groups according to the nomenclature of \cite{carat}.
Table \ref{tab:Number-Totally-Zone-Cont}
and \ref{tab:Table-Irreducible-Domains} 
relate our classification to notions in the theory of {\em contraction types}
as developed in \cite{engel-2000}. In the following we provide some
background information (see also \cite{SumParallZonotop}).

\subsection{Fundamental faces and irreducible cones}

For a given secondary cone $\SC$ with generating rays $R_1, \dots, R_k$ we define
the {\em fundamental face} $F(\SC)$ to be the smallest face of $\SC$ that
contains all the generators $R_i$ of rank greater than~$1$.
The face $F(\SC)$ may be reduced to zero in which case $\SC$ is generated
by rank-$1$ matrices only.
From \cite{er-1994} we know that the number of generators is equal to
the dimension of the secondary cone in this case and that this case is equivalent to the Dirichlet-Voronoi polytope being
a zonotope and to the Delaunay subdivision being the connected region
of a hyperplane arrangement.
Up to $\GL_5(\Z)$-equivalence,
we found $81$ secondary cones of this kind, 
corresponding to different zonotopes in dimension~$5$.

If $F(\SC)$ is nontrivial (non-zero) then the structure of the
secondary cone is more complex.
For a secondary cone $\SC$ we have a decomposition of the form
\begin{equation*}
\SC = F(\SC) + \sum_{i=1}^h \R_+ p(v_i) 
,
\end{equation*}
with $p(v_i)=v_i v_i^t$ the rank-$1$ matrix (form) associated to a vector~$v_i$.
Our computations show that we have $\dim \SC = \dim F(\SC) + h$ which means that $\SC$ is obtained
by a sequence of $h$ pyramid constructions over $F(\SC)$.
By a pyramid construction we mean an extension to a 
higher dimensional secondary cone by adding a rank-$1$ generating ray.

If $F(\SC)$ does not contain any positive definite matrices 
(and hence lies in the boundary of ${\mathcal S}^5_{>0}$),
then in dimension~$5$
there is only one possibility: $F(\SC)$ has only one 
extreme ray that corresponds
to the $\mathsf{D}_4$ root lattice, which we denote by $F_{\mathsf{D}_4}$.
Up to $\GL_5(\Z)$-equivalence, we found $424$ different combinatorial
types of secondary cones of the form $F_{\mathsf{D}_4} + \sum_{i=1}^h \R_+ p(v_i)$.
Note that $F_{\mathsf{D}_4}$ itself is not a secondary cone, since it does
not contain any positive definite forms.
By our computation, all such cones have their dimension equal to
their number of generators. 

The fundamental cones $F(\SC)$ may themselves contain rank $1$-forms.
For example, there exist two secondary cones of dimension $3$ with $4$
generators each, $3$ of rank $4$
and one of rank $1$ (see Section 5 of \cite{SmoothnessRegularity}).
If $F(\SC)$ contains only forms of rank higher than $1$ then according to the terminology
of \cite{engel-2000} it is {\em totally zone contracted}. If a secondary cone satisfies
$\SC = F(\SC)$ then it is called {\em irreducible}. Table~\ref{tab:Number-Totally-Zone-Cont}
and~\ref{tab:Table-Irreducible-Domains} give key information on
irreducible secondary cones we found.

\subsection{Contraction types}

In \cite{engel-2000} the notion of a contraction type is introduced. This notion is distinct
from secondary cones and gives a further refinement of them. That is, if we have a secondary
cone $\SC$ that is irreducible but not totally zone-contracted and has rank-$1$ forms $p_1$, \dots, $p_m$,
then we can decompose it into a number of 
{\em contraction cones} (also called {\em contraction domains}) 
$\SC_i + \sum_{j=1}^m \R_+ p_j$ 
with $\SC_i$ a totally zone-contracted secondary cone.
For example the $3$-dimensional cone $\SC$ with symbol $L_1^2L_3 p_1$ 
in Table~\ref{tab:Table-Irreducible-Domains} is a cone over a square (combinatorially) with vertices
corresponding to $p_1$, $L_1$, $L_3$ and $L_1$. We can decompose it into
two isomorphic $3$-dimensional cones (over triangles) of the form $L_1L_3 + \R_+ p_1$ and
one $2$-dimensional cone of the form $L_3 + \R_+ p_1$.

For other cones the decomposition can be more complicated. 
Given an irreducible secondary cone $\SC$, 
let $R_1$ be the cone 
of its extreme rays of rank~$1$. 
We define ${\mathcal S}$ to be the set of all 
totally zone contracted irreducible cones 
whose rays are also rays of $\SC$
(of rank greater than $1$).
Then our computation shows that $\SC$
can be decomposed into contraction cones $S + R_1$ with $S\in {\mathcal S}$.

The decomposition of an irreducible secondary cone~$\SC$
into contraction cones, induces a decomposition of any secondary 
cone obtained by adding rank-$1$ forms.
Overall, we thus obtain a decomposition into contraction cones that is finer than the
decomposition by secondary cones. For secondary cones $\SC$ whose fundamental face $F(\SC)$
is totally zone-contracted there is no difference.
But for other irreducible secondary cones the contraction
types form a strictly finer decomposition. The total number of contraction types that we obtain
is~$181394$.
The number of contraction cones by their dimension is given in
Table \ref{tab:number-secondary-cones-contraction-cone-bydim}.
In Table~\ref{tab:Table-Irreducible-Domains} we give for each
irreducible secondary cone~$D$
the number of types of 
contraction cones contained in $D + \sum_k \R_+ p(v_k)$.
We note that in \cite{engel-2000} the number of contraction cones is reported
to be $179372$. This discrepancy is most likely due to
the different notion of equivalence via ``subordination schemes'' used there.

\subsection{Euler Poincar\'e characteristic check}

Another key check of the correctness of our enumeration is to use
the Euler Poincar\'e characteristic.
We have the formula 
\begin{equation*}
\sum_{F} (-1)^{\dim(F)} \frac{1}{\vert\Stab(F)\vert} = 0
\end{equation*}
where the sum is over the representatives of cones with respect 
to the action of $\GL_n(\Z)$.
This kind of formula comes from the Euler Poincar\'e characteristic of discrete
groups, i.e. $\chi(\GL_n(\Z))=0$ for $n\geq 3$.
See \cite{BrownCohobook,CohomologyImaginary} for more details.

Both, our enumeration of secondary cones and our enumeration of contraction
cones satisfy this condition, 
which is yet another strong indication of the correctness of our enumeration.
For example for the secondary cones, if we regroup the cones by their dimension,
this gives us the following non-trivial identity:
\begin{equation*}
\begin{array}{c}
- \frac{293}{5760} + \frac{7463}{5760} - \frac{939}{64} + \frac{56927}{576} - \frac{5146751}{11520}\\
+ \frac{8329297}{5760} - \frac{3341911}{960} + \frac{1630783}{256} - \frac{10308319}{1152} +
\frac{13879537}{1440}\\
- \frac{1414553}{180} + \frac{1356727}{288} - \frac{565595}{288} + \frac{48907}{96} - \frac{8923}{144} = 0.
\end{array}
\end{equation*}
This kind of mass formula provides a highly non-trivial check 
of the correctness of an enumeration
as any error on a single entry or on a single stabilizer would turn the formula wrong.

\section{Acknowledgements}

The authors like to thank Peter Engel for several
helpful communications regarding the types
classified in \cite{engel-2000}, according to ``subordination schemes''.
Mathieu Dutour Sikiri\'c and Achill Sch\"urmann were partially supported 
by the Erwin-Schr\"odinger-Institute (ESI) during a stay in fall 2014 for
the program on Minimal Energy Point Sets, Lattices and Designs.
Their research was also supported by the Humboldt foundation
and DFG-grant SCHU-1503/6-1.


\providecommand{\bysame}{\leavevmode\hbox to3em{\hrulefill}\thinspace}
\providecommand{\MR}{\relax\ifhmode\unskip\space\fi MR }
\providecommand{\MRhref}[2]{%
  \href{http://www.ams.org/mathscinet-getitem?mr=#1}{#2}
}
\providecommand{\href}[2]{#2}



\newcommand{\etalchar}[1]{$^{#1}$}


\newpage

\begin{table}
\caption{Number of $\GL_5(\Z)$-inequivalent secondary cones and contraction cones in~${\mathcal S}^5_{>0}$ by their dimension.}
\label{tab:number-secondary-cones-contraction-cone-bydim}
\begin{center}
\begin{tabular}{|c|c|c||c|c|c|}
  \cline{1-6}
  n & nr. sec. c. & nr. cont. c. & n & nr. sec. c. & nr. cont. c.\\
  \cline{1-6}
  1 & 7 & 7 & 9 & 21132 & 33085\\
  2 & 37 & 39 & 10 & 22221 & 37601\\
  3 & 146 & 161 & 11 & 18033 & 32821\\
  4 & 535 & 613 & 12 & 10886 & 21292\\
  5 & 1681 & 2021 & 13 & 4713 & 9709\\
  6 & 4366 & 5543 & 14 & 1318 & 2787\\
  7 & 9255 & 12512 & 15 & 222 & 397\\
  8 & 15692 & 22806 &  &  &\\
  \cline{1-6}
\end{tabular}
\end{center}

\end{table}

\begin{table}
\caption{Number of $\GL_5(\Z)$-inequivalent secondary cones
  in~${\mathcal S}^5_{>0}$ by number of rank-$k$ generating rays.
  In line $i$, the rank-$k$-column, $k=1,4,5$, contains the number of
  secondary cones which have $i$ generating rays of rank $k$.
 (There exist no generating rays for $k=2,3$.)}
\label{tab:number-secondary-cones-byrank}
\begin{tabular}{lrrr}
\parbox[b]{3cm}{\# Generating rays \\ $\,$(of particular rank)}
& rank-$1$ & rank-$4$ & rank-$5$ \\
\hline 
$0$ & $82$ & $51900$ & $1572$ \\
$1$ & $410$ & $35316$ & $15421$ \\
$2$ & $1658$ & $21574$ & $32939$ \\
$3$ & $5029$ &  $1354$ & $26811$ \\
$4$ & $11301$ & $0$ & $19302$ \\
$5$ & $18923$ &  $100$ & $6841$ \\
$6$ & $23802$ & $0$ & $3662$ \\
$7$ & $22411$ & $0$ & $2150$ \\
$8$ & $15528$ & $0$ & $950$ \\
$9$ & $7744$ & $0$ & $285$ \\
$10$ & $2699$ & $0$ & $170$ \\
$11$ & $548$ & $0$ & $38$ \\
$12$ & $97$ & $0$ & $76$ \\
$13$ & $9$ & $0$ & $0$ \\
$14$ & $2$ & $0$ & $0$ \\
$15$ & $1$ & $0$ & $9$ \\
$16$ & $0$ & $0$ & $18$ \\
\hline 
\end{tabular}
\end{table}

\onecolumn

\begin{table}
  \caption{Number of primitive and all combinatorial types of Delaunay
    subdivisions, respectively corresponding $\GL_n(\Z)$-inequivalent secondary cones.}
  \label{tab:NumberLtypeEnumeration}
  \begin{center}
    \begin{tabular}{|c|c|c|}
      \cline{1-3}
      n & Primitive types & All combinatorial types\\
      \cline{1-3}
      2 & $1$                                       & $2$\\
      3 & $1$ \cite{fedorov-1885}                   & $5$ \cite{fedorov-1885}\\
      4 & $3$ \cite{voronoi-1908}                   & $52$ \cite{delone-1929,stogrin-1973}\\
      5 & $222$ \cite{br-1973}                      & $110244$ \\
        & \cite{rb-1976,eg-2002}                    &\\
      6 & $\geq 567.613.632$ \cite{be-2013}         & \\
      \cline{1-3}
    \end{tabular}
  \end{center}
\end{table}

\begin{table}
\centering
\caption{Number of secondary cones according to dimension (at most $15$) and number of generators (at most $26$).}
\label{tab:number-domain-by-dim-and-gen}
\begin{center}
  \begin{tabular}{|c|ccccccccccccccc|}
    \cline{1-16}
    & 1 & 2 & 3 & 4 & 5 & 6 & 7 & 8 & 9 & 10 & 11 & 12 & 13 & 14 & 15\\
    \cline{1-16}
    1 &7 & & & & & & & & & & & & & &\\
    2 & &37 & & & & & & & & & & & & &\\
    3 & & &144 & & & & & & & & & & & &\\
    4 & & &2 &517 & & & & & & & & & & &\\
    5 & & & &17 &1595 & & & & & & & & & &\\
    6 & & & & &81 &4041 & & & & & & & & &\\
    7 & & & & &1 &301 &8266 & & & & & & & &\\
    8 & & & &1 & &12 &887 &13354 & & & & & & &\\
    9 & & & & &3 & &62 &2007 &16862 & & & & & &\\
    10 & & & & &1 &11 &2 &222 &3461 &16358 & & & & &\\
    11 & & & & & &1 &36 &13 &557 &4443 &11989 & & & &\\
    12 & & & & & & &2 &89 &50 &944 &4259 &6395 & & &\\
    13 & & & & & & & &7 &182 &122 &1103 &2945 &2346 & &\\
    14 & & & & & & & & &19 &305 &181 &857 &1449 &526 &\\
    15 & & & & & & & & & &43 &403 &173 &430 &456 &62\\
    16 & & & & & & & & &1 & &80 &390 &102 &120 &84\\
    17 & & & & & & & & & &5 & &92 &274 &35 &13\\
    18 & & & & & & & & & & &15 & &72 &122 &5\\
    19 & & & & & & & & & & & &30 & &29 &33\\
    20 & & & & & & & & & & & & &34 & &13\\
    21 & & & & & & & & & &1 & & & &23 &\\
    22 & & & & & & & & & & &3 & & & &6\\
    23 & & & & & & & & & & & &4 & & &\\
    24 & & & & & & & & & & & & &6 & &\\
    25 & & & & & & & & & & & & & &7 &\\
    26 & & & & & & & & & & & & & & &6\\
    \cline{1-16}
  \end{tabular}
\end{center}
\end{table}

\begin{table}
\caption{Number of $\GL_5(\Z)$-inequivalent secondary cones in~${\mathcal S}^5_{>0}$ 
which are not extendable to a higher dimensional secondary cone by adding a
rank-$1$ generating ray.}
\label{tab:number-secondary-cones-not-extendable}
\begin{tabular}{lrrrrrr}
Dimension & $10$ & $11$ & $12$ & $13$ & $14$ & $15$ \\ 
\hline
\# Secondary cones & $1$ & $12$ & $40$ & $142$ & $266$ & $222$\\
\end{tabular}
\end{table}

\begin{table}
\caption{Frequency of occurence of Bravais groups. ``name'' is the
  standard name from the GAP package \cite{carat}. ``order'' is the
  size of the point group of corresponding lattices. ``frequency'' is
  the number of secondary cones that are symmetric with respect to the
  group.}
\label{tab:freq-bravais-group}
\begin{center}
  \begin{tabular}{|c|c|c||c|c|c||c|c|c|}
    \cline{1-9}
    name & order & frequency & name & order & frequency & name & order & frequency\\
    \cline{1-9}
    1,1,1,1,1 :1 & 2 & 105301 & 1,1;1;1;1 :17 & 16 & 5 & 4-1;1 :2 & 768 & 1\\
    1,1,1,1;1 :2 & 4 & 4155 & 3;1,1 :2 & 96 & 4 & 4-1;1 :3 & 2304 & 1\\
    1,1,1;1;1 :6 & 8 & 159 & 3;1,1 :5 & 96 & 4 & 5-1 :3 & 3840 & 1\\
    2-2;1,1,1 :2 & 12 & 137 & 2-1;1,1;1 :6 & 32 & 4 & 5-2 :3 & 1440 & 1\\
    1,1,1;1,1 :2 & 4 & 112 & 1;1;1;1;1 :8 & 32 & 4 & 3;1;1 :4 & 192 & 1\\
    1,1,1;1;1 :4 & 8 & 90 & 1,1,1;1,1 :1 & 4 & 3 & 4-1;1 :4 & 768 & 1\\
    1,1,1;1;1 :5 & 8 & 39 & 1,1,1;1;1 :1 & 8 & 3 & 2-2;2-2;1 :5 & 72 & 1\\
    1,1,1,1;1 :1 & 4 & 34 & 2-2;2-2;1 :3 & 72 & 3 & 2-1;1;1;1 :6 & 64 & 1\\
    2-1;1,1,1 :2 & 16 & 31 & 1,1;1;1;1 :10 & 16 & 3 & 2-1;1;1;1 :7 & 64 & 1\\
    2-2;1,1;1 :6 & 24 & 31 & 4-3;1 :3 & 240 & 2 & 2-2;1;1;1 :7 & 48 & 1\\
    1,1;1;1;1 :15 & 16 & 20 & 2-2;1,1;1 :4 & 24 & 2 & 3;1;1 :7 & 192 & 1\\
    1,1;1,1;1 :3 & 8 & 14 & 1;1;1;1;1 :5 & 32 & 2 & 2-1;1;1;1 :8 & 64 & 1\\
    1,1;1;1;1 :13 & 16 & 12 & 2-2;1,1;1 :5 & 24 & 2 & 2-1;1;1;1 :11 & 64 & 1\\
    3;1,1 :3 & 48 & 10 & 3;1;1 :12 & 192 & 2 & 1;1;1;1;1 :12 & 32 & 1\\
    1,1;1;1;1 :6 & 16 & 8 & 1;1;1;1;1 :13 & 32 & 2 & 2-1;1;1;1 :12 & 64 & 1\\
    3;1;1 :8 & 96 & 7 & 1,1;1,1;1 :1 & 8 & 1 & 1;1;1;1;1 :15 & 32 & 1\\
    1,1,1;1;1 :2 & 8 & 6 & 1,1;1;1;1 :1 & 16 & 1 & 1;1;1;1;1 :16 & 32 & 1\\
    2-1;1,1;1 :4 & 32 & 6 & 1;1;1;1;1 :1 & 32 & 1 &  & & \\
    1,1;1,1;1 :6 & 8 & 6 & 3;1;1 :2 & 192 & 1 &  & & \\
    \cline{1-9}
  \end{tabular}
\end{center}

\end{table}

\begin{table}
\caption{Information about the $82$ totally zone-contracted secondary cones. ``dim'' is the
  dimension of the secondary cone~$SC$, ``generator'' gives the type of
  the extreme rays, ``symbol'' gives the number of facets and vertices
  of the corresponding Dirichlet-Voronoi polytopes
  and ``nb sec. c.'' gives the number of secondary cones having
  $SC$ as their fundamental face.
}
\label{tab:Number-Totally-Zone-Cont}
\begin{scriptsize}  
\begin{center}
\begin{tabular}{|c|c|c|c||c|c|c|c||c|c|c|c|}
\cline{1-12}
dim & generator & symbol & nb sec. c. & dim & generator & symbol & nb sec. c. & dim & generator & symbol & nb sec. c.\\
\cline{1-12}
1 & $L_{1}$ & 40,42 & 450 & 3 & $L_{1}\mathsf{D}_{4}^{2}$ & 48,242 & 2738 & 4 & $L_{2}^{2}\mathsf{D}_{4}^{3}$ & 42,204 & 665\\
1 & $L_{2}$ & 42,96 & 777 & 3 & $L_{2}\mathsf{D}_{4}^{2}$ & 42,168 & 2047 & 4 & $L_{1}L_{2}\mathsf{D}_{4}^{2}$ & 48,282 & 3988\\
1 & $L_{3}$ & 48,180 & 670 & 3 & $L_{3}\mathsf{D}_{4}^{2}$ & 52,344 & 1344 & 4 & $L_{1}L_{3}\mathsf{D}_{4}^{2}$ & 52,352 & 2272\\
1 & $L_{4}$ & 50,192 & 112 & 3 & $L_{7}\mathsf{D}_{4}^{2}$ & 56,462 & 484 & 4 & $L_{2}L_{3}\mathsf{D}_{4}^{2}$ & 52,384 & 1074\\
1 & $L_{5}$ & 50,282 & 352 & 3 & $L_{1}L_{2}\mathsf{D}_{4}$ & 48,242 & 5029 & 4 & $L_{3}L_{7}\mathsf{D}_{4}^{2}$ & 56,470 & 1160\\
1 & $L_{6}$ & 54,342 & 324 & 3 & $L_{1}L_{3}\mathsf{D}_{4}$ & 48,254 & 2436 & 4 & $L_{1}L_{2}L_{3}\mathsf{D}_{4}$ & 52,354 & 4100\\
1 & $L_{7}$ & 54,366 & 220 & 3 & $L_{1}L_{5}\mathsf{D}_{4}$ & 50,328 & 650 & 4 & $L_{1}L_{2}L_{5}\mathsf{D}_{4}$ & 54,418 & 1256\\
2 & $\mathsf{D}_{4}^{2}$ & 42,132 & 1067 & 3 & $L_{2}L_{3}\mathsf{D}_{4}$ & 52,346 & 2344 & 4 & $L_{1}L_{3}^{2}L_{7}\mathsf{D}_{4}$ & 54,418 & 1088\\
2 & $L_{1}\mathsf{D}_{4}$ & 40,122 & 1814 & 3 & $L_{2}L_{5}\mathsf{D}_{4}$ & 54,402 & 650 & 4 & $L_{1}L_{3}L_{5}\mathsf{D}_{4}$ & 50,342 & 696\\
2 & $L_{2}\mathsf{D}_{4}$ & 42,132 & 1825 & 3 & $L_{3}L_{5}\mathsf{D}_{4}$ & 50,334 & 553 & 4 & $L_{2}L_{3}L_{5}\mathsf{D}_{4}$ & 54,424 & 1092\\
2 & $L_{3}\mathsf{D}_{4}$ & 48,246 & 1428 & 3 & $L_{3}L_{7}\mathsf{D}_{4}$ & 54,410 & 1160 & 4 & $L_{1}L_{2}L_{3}L_{5}$ & 54,406 & 1392\\
2 & $L_{5}\mathsf{D}_{4}$ & 50,312 & 352 & 3 & $L_{1}L_{2}L_{3}$ & 52,316 & 2773 & 4 & $L_{1}L_{2}L_{3}L_{6}$ & 54,428 & 856\\
2 & $L_{7}\mathsf{D}_{4}$ & 54,402 & 484 & 3 & $L_{1}L_{2}L_{5}$ & 54,392 & 1256 & 4 & $L_{1}L_{2}L_{5}L_{6}$ & 54,438 & 928\\
2 & $L_{1}L_{2}$ & 48,202 & 2385 & 3 & $L_{1}L_{2}L_{6}$ & 54,400 & 758 & 4 & $L_{1}L_{3}L_{4}L_{5}$ & 50,360 & 696\\
2 & $L_{1}L_{3}$ & 48,188 & 1058 & 3 & $L_{1}L_{3}^{2}L_{7}$ & 54,382 & 456 & 4 & $L_{1}L_{3}L_{4}L_{6}$ & 54,416 & 786\\
2 & $L_{1}L_{4}$ & 50,232 & 333 & 3 & $L_{1}L_{3}L_{4}$ & 50,288 & 516 & 4 & $L_{1}L_{3}L_{5}L_{6}$ & 54,418 & 800\\
2 & $L_{1}L_{5}$ & 50,298 & 650 & 3 & $L_{1}L_{3}L_{5}$ & 50,312 & 696 & 4 & $L_{1}L_{4}L_{5}L_{6}$ & 54,426 & 928\\
2 & $L_{1}L_{6}$ & 54,366 & 758 & 3 & $L_{1}L_{3}L_{6}$ & 54,394 & 856 & 4 & $L_{2}L_{3}L_{5}L_{6}$ & 54,444 & 628\\
2 & $L_{2}L_{3}$ & 52,308 & 1638 & 3 & $L_{1}L_{4}L_{5}$ & 50,346 & 630 & 4 & $L_{3}L_{4}L_{5}L_{6}$ & 54,432 & 628\\
2 & $L_{2}L_{5}$ & 54,376 & 650 & 3 & $L_{1}L_{4}L_{6}$ & 54,388 & 734 & 5 & $L_{2}^{5}\mathsf{D}_{4}^{5}$ & 42,240 & 100\\
2 & $L_{2}L_{6}$ & 54,376 & 324 & 3 & $L_{1}L_{5}L_{6}$ & 54,404 & 928 & 5 & $L_{1}L_{2}^{2}\mathsf{D}_{4}^{3}$ & 48,322 & 689\\
2 & $L_{3}L_{4}$ & 50,280 & 318 & 3 & $L_{2}L_{3}L_{5}$ & 54,398 & 1092 & 5 & $L_{1}L_{2}L_{3}\mathsf{D}_{4}^{2}$ & 52,392 & 1815\\
2 & $L_{3}L_{5}$ & 50,304 & 553 & 3 & $L_{2}L_{3}L_{6}$ & 54,420 & 582 & 5 & $L_{1}L_{3}^{2}L_{7}\mathsf{D}_{4}^{2}$ & 56,478 & 1088\\
2 & $L_{3}L_{6}$ & 54,386 & 582 & 3 & $L_{2}L_{5}L_{6}$ & 54,422 & 553 & 5 & $L_{1}L_{2}L_{3}L_{5}\mathsf{D}_{4}$ & 54,432 & 1392\\
2 & $L_{3}L_{7}$ & 54,374 & 490 & 3 & $L_{3}L_{4}L_{5}$ & 50,352 & 553 & 5 & $L_{1}L_{2}L_{3}L_{5}L_{6}$ & 54,452 & 800\\
2 & $L_{4}L_{5}$ & 50,330 & 348 & 3 & $L_{3}L_{4}L_{6}$ & 54,408 & 531 & 5 & $L_{1}L_{3}L_{4}L_{5}L_{6}$ & 54,440 & 800\\
2 & $L_{4}L_{6}$ & 54,364 & 318 & 3 & $L_{3}L_{5}L_{6}$ & 54,410 & 628 &  & & & \\
2 & $L_{5}L_{6}$ & 54,388 & 553 & 3 & $L_{4}L_{5}L_{6}$ & 54,410 & 553 &  & & & \\
\cline{1-12}
\end{tabular}
\end{center}
\end{scriptsize}
\end{table}

\begin{table}
\caption{Information about the $125$ inequivalent irreducible secondary cones,
  which are not totally zone-contracted. 
Same convention as in Table~\ref{tab:Number-Totally-Zone-Cont}; 
in addition $p_1$ denotes an extreme ray of rank $1$ and ``nb cont. d.'' is the number of contraction cones corresponding to this irreducible component.}
\label{tab:Table-Irreducible-Domains}
\begin{scriptsize}  
\begin{center}
\begin{tabular}{|c|c|c|c|c||c|c|c|c|c|}
\cline{1-10}
dim & generator & symbol & nb sec. c. & nb cont. c. & dim & generator & symbol & nb sec. c. & nb cont. c.\\
\cline{1-10}
3 & $L_{1}^{2}L_{3}p_1$ & 48,196 & 566 & 2047 & 7 & $L_{3}L_{4}L_{6}^{2}p_1^{4}$ & 58,536 & 27 & 72\\
4 & $L_{1}L_{3}L_{5}p_1^{2}$ & 50,320 & 205 & 3988 & 7 & $L_{1}^{3}L_{3}^{3}L_{5}\mathsf{D}_{4}p_1^{3}$ & 50,378 & 73 & 3030\\
4 & $L_{1}^{2}L_{3}\mathsf{D}_{4}p_1$ & 48,262 & 1240 & 1074 & 7 & $L_{1}^{3}L_{2}L_{3}^{3}L_{5}p_1^{3}$ & 54,442 & 134 & 639\\
4 & $L_{1}^{3}L_{3}^{3}L_{7}p_1$ & 54,390 & 174 & 665 & 7 & $L_{1}^{3}L_{3}^{3}L_{4}L_{5}p_1^{3}$ & 50,396 & 73 & 1274\\
4 & $L_{1}^{2}L_{2}L_{3}p_1$ & 52,324 & 1423 & 1092 & 7 & $L_{1}^{3}L_{3}^{3}L_{5}L_{6}p_1^{3}$ & 54,454 & 33 & 820\\
4 & $L_{1}^{2}L_{3}L_{4}p_1$ & 50,296 & 274 & 1256 & 7 & $L_{1}^{2}L_{2}L_{3}L_{5}L_{6}p_1^{3}$ & 54,474 & 164 & 605\\
4 & $L_{1}^{2}L_{3}L_{5}p_1$ & 50,320 & 205 & 615 & 7 & $L_{1}^{2}L_{3}^{2}L_{5}^{2}L_{6}p_1^{3}$ & 54,464 & 74 & 1000\\
4 & $L_{1}^{2}L_{3}L_{6}p_1$ & 54,402 & 358 & 4100 & 7 & $L_{1}^{2}L_{3}^{2}L_{5}^{2}L_{6}p_1^{3}$ & 54,464 & 148 & 740\\
4 & $L_{1}L_{3}^{2}L_{5}p_1$ & 50,326 & 182 & 3503 & 7 & $L_{1}^{2}L_{3}L_{4}L_{5}L_{6}p_1^{3}$ & 54,462 & 150 & 207\\
4 & $L_{3}L_{5}^{2}L_{6}p_1$ & 54,434 & 203 & 3999 & 7 & $L_{1}L_{2}L_{3}L_{5}^{2}L_{6}p_1^{3}$ & 54,484 & 121 & 814\\
5 & $L_{1}L_{5}L_{6}p_1^{3}$ & 54,412 & 97 & 615 & 7 & $L_{1}L_{3}^{3}L_{5}^{3}L_{6}p_1^{3}$ & 54,480 & 22 & 261\\
5 & $L_{1}L_{3}L_{5}\mathsf{D}_{4}p_1^{2}$ & 50,350 & 205 & 1188 & 7 & $L_{1}L_{3}L_{4}L_{5}^{2}L_{6}p_1^{3}$ & 54,472 & 121 & 1036\\
5 & $L_{1}^{2}L_{3}^{2}L_{5}p_1^{2}$ & 50,334 & 298 & 5895 & 7 & $L_{1}^{2}L_{2}L_{3}^{2}L_{5}\mathsf{D}_{4}p_1^{2}$ & 54,454 & 606 & 153\\
5 & $L_{1}L_{2}L_{3}L_{5}p_1^{2}$ & 54,414 & 396 & 492 & 7 & $L_{1}^{2}L_{2}L_{3}^{2}L_{5}L_{6}p_1^{2}$ & 54,474 & 200 & 750\\
5 & $L_{1}L_{3}L_{4}L_{5}p_1^{2}$ & 50,368 & 197 & 492 & 7 & $L_{1}^{2}L_{2}L_{3}L_{5}^{2}L_{6}p_1^{2}$ & 54,484 & 34 & 605\\
5 & $L_{1}L_{3}L_{5}L_{6}p_1^{2}$ & 54,426 & 164 & 689 & 7 & $L_{1}^{2}L_{3}^{2}L_{4}L_{5}L_{6}p_1^{2}$ & 54,462 & 200 & 1000\\
5 & $L_{1}L_{3}L_{5}L_{6}p_1^{2}$ & 54,432 & 164 & 1815 & 7 & $L_{1}^{2}L_{3}L_{4}L_{5}^{2}L_{6}p_1^{2}$ & 54,472 & 34 & 740\\
5 & $L_{1}^{2}L_{3}\mathsf{D}_{4}^{2}p_1$ & 52,360 & 1168 & 3279 & 7 & $L_{1}L_{2}L_{3}^{2}L_{5}^{2}L_{6}p_1^{2}$ & 54,490 & 148 & 207\\
5 & $L_{1}^{3}L_{3}^{3}L_{7}\mathsf{D}_{4}p_1$ & 54,426 & 396 & 100 & 7 & $L_{1}L_{3}^{2}L_{4}L_{5}^{2}L_{6}p_1^{2}$ & 54,478 & 148 & 639\\
5 & $L_{1}^{2}L_{2}L_{3}\mathsf{D}_{4}p_1$ & 52,362 & 2060 & 1392 & 8 & $L_{1}^{2}L_{3}^{2}L_{5}^{2}L_{6}p_1^{5}$ & 54,478 & 34 & 320\\
5 & $L_{1}^{2}L_{3}L_{5}\mathsf{D}_{4}p_1$ & 50,350 & 205 & 553 & 8 & $L_{1}^{3}L_{3}^{3}L_{5}^{2}L_{6}p_1^{4}$ & 54,478 & 47 & 1274\\
5 & $L_{1}L_{3}^{2}L_{5}\mathsf{D}_{4}p_1$ & 50,356 & 182 & 1092 & 8 & $L_{1}^{2}L_{3}^{3}L_{5}^{3}L_{6}p_1^{4}$ & 54,488 & 38 & 337\\
5 & $L_{1}^{2}L_{2}L_{3}L_{5}p_1$ & 54,414 & 396 & 958 & 8 & $L_{1}^{2}L_{3}^{3}L_{5}^{2}L_{6}p_1^{4}$ & 54,478 & 43 & 814\\
5 & $L_{1}^{2}L_{2}L_{3}L_{6}p_1$ & 54,436 & 358 & 480 & 8 & $L_{1}L_{2}L_{3}^{2}L_{5}L_{6}p_1^{4}$ & 54,488 & 64 & 487\\
5 & $L_{1}^{2}L_{3}L_{4}L_{5}p_1$ & 50,368 & 205 & 1490 & 8 & $L_{1}L_{3}^{2}L_{4}L_{5}L_{6}p_1^{4}$ & 54,476 & 57 & 285\\
5 & $L_{1}^{2}L_{3}L_{4}L_{6}p_1$ & 54,424 & 327 & 990 & 8 & $L_{1}L_{3}L_{4}L_{6}^{2}p_1^{4}$ & 58,544 & 28 & 77\\
5 & $L_{1}^{2}L_{3}L_{5}L_{6}p_1$ & 54,426 & 228 & 291 & 8 & $L_{1}^{3}L_{2}L_{3}^{3}L_{5}\mathsf{D}_{4}p_1^{3}$ & 54,468 & 134 & 261\\
5 & $L_{1}L_{2}L_{3}^{2}L_{5}p_1$ & 54,420 & 352 & 546 & 8 & $L_{1}^{3}L_{2}L_{3}^{3}L_{5}L_{6}p_1^{3}$ & 54,488 & 33 & 1036\\
5 & $L_{1}L_{3}^{2}L_{4}L_{5}p_1$ & 50,374 & 182 & 800 & 8 & $L_{1}^{3}L_{3}^{3}L_{4}L_{5}L_{6}p_1^{3}$ & 54,476 & 33 & 753\\
5 & $L_{1}L_{3}^{2}L_{5}L_{6}p_1$ & 54,432 & 128 & 628 & 8 & $L_{1}^{2}L_{2}L_{3}^{2}L_{5}^{2}L_{6}p_1^{3}$ & 54,498 & 74 & 153\\
5 & $L_{1}L_{3}L_{5}^{2}L_{6}p_1$ & 54,442 & 178 & 328 & 8 & $L_{1}^{2}L_{2}L_{3}^{2}L_{5}^{2}L_{6}p_1^{3}$ & 54,498 & 148 & 575\\
5 & $L_{2}L_{3}L_{5}^{2}L_{6}p_1$ & 54,468 & 203 & 474 & 8 & $L_{1}^{2}L_{3}^{2}L_{4}L_{5}^{2}L_{6}p_1^{3}$ & 54,486 & 74 & 814\\
5 & $L_{3}L_{4}L_{5}^{2}L_{6}p_1$ & 54,456 & 203 & 591 & 8 & $L_{1}^{2}L_{3}^{2}L_{4}L_{5}^{2}L_{6}p_1^{3}$ & 54,486 & 148 & 261\\
6 & $L_{3}^{2}L_{6}p_1^{4}$ & 54,430 & 34 & 92 & 8 & $L_{1}L_{2}L_{3}^{3}L_{5}^{3}L_{6}p_1^{3}$ & 54,514 & 22 & 1036\\
6 & $L_{1}^{3}L_{3}^{3}L_{5}p_1^{3}$ & 50,348 & 73 & 1188 & 8 & $L_{1}L_{3}^{3}L_{4}L_{5}^{3}L_{6}p_1^{3}$ & 54,502 & 22 & 153\\
6 & $L_{1}^{2}L_{3}L_{5}L_{6}p_1^{3}$ & 54,440 & 164 & 492 & 9 & $L_{1}L_{4}L_{5}L_{6}p_1^{6}$ & 54,502 & 16 & 337\\
6 & $L_{1}L_{2}L_{5}L_{6}p_1^{3}$ & 54,446 & 97 & 492 & 9 & $L_{1}^{3}L_{3}^{4}L_{5}^{3}L_{6}p_1^{5}$ & 54,502 & 38 & 487\\
6 & $L_{1}L_{3}L_{5}^{2}L_{6}p_1^{3}$ & 54,450 & 121 & 2619 & 9 & $L_{1}^{2}L_{2}L_{3}^{2}L_{5}^{2}L_{6}p_1^{5}$ & 54,512 & 34 & 48\\
6 & $L_{1}L_{4}L_{5}L_{6}p_1^{3}$ & 54,434 & 93 & 1092 & 9 & $L_{1}^{2}L_{3}^{2}L_{4}L_{5}^{2}L_{6}p_1^{5}$ & 54,500 & 30 & 753\\
6 & $L_{1}^{2}L_{3}^{2}L_{5}\mathsf{D}_{4}p_1^{2}$ & 50,364 & 298 & 958 & 9 & $L_{1}^{2}L_{3}L_{4}L_{6}^{2}p_1^{5}$ & 58,552 & 11 & 575\\
6 & $L_{1}L_{2}L_{3}L_{5}\mathsf{D}_{4}p_1^{2}$ & 54,440 & 396 & 1490 & 9 & $L_{1}^{3}L_{2}L_{3}^{3}L_{5}^{2}L_{6}p_1^{4}$ & 54,512 & 47 & 905\\
6 & $L_{1}^{2}L_{2}L_{3}^{2}L_{5}p_1^{2}$ & 54,428 & 606 & 3030 & 9 & $L_{1}^{3}L_{3}^{3}L_{4}L_{5}^{2}L_{6}p_1^{4}$ & 54,500 & 47 & 300\\
6 & $L_{1}^{2}L_{3}^{2}L_{4}L_{5}p_1^{2}$ & 50,382 & 298 & 639 & 9 & $L_{1}^{2}L_{2}L_{3}^{3}L_{5}^{3}L_{6}p_1^{4}$ & 54,522 & 38 & 68\\
6 & $L_{1}^{2}L_{3}^{2}L_{5}L_{6}p_1^{2}$ & 54,440 & 200 & 291 & 9 & $L_{1}^{2}L_{2}L_{3}^{3}L_{5}^{2}L_{6}p_1^{4}$ & 54,512 & 43 & 487\\
6 & $L_{1}^{2}L_{3}L_{5}^{2}L_{6}p_1^{2}$ & 54,450 & 34 & 820 & 9 & $L_{1}^{2}L_{3}^{3}L_{4}L_{5}^{3}L_{6}p_1^{4}$ & 54,510 & 38 & 753\\
6 & $L_{1}L_{2}L_{3}L_{5}L_{6}p_1^{2}$ & 54,460 & 164 & 605 & 9 & $L_{1}^{2}L_{3}^{3}L_{4}L_{5}^{2}L_{6}p_1^{4}$ & 54,500 & 43 & 575\\
6 & $L_{1}L_{2}L_{3}L_{5}L_{6}p_1^{2}$ & 54,466 & 164 & 628 & 10 & $L_{3}L_{4}L_{6}p_1^{8}$ & 54,452 & 6 & 18\\
6 & $L_{1}L_{3}^{2}L_{5}^{2}L_{6}p_1^{2}$ & 54,456 & 148 & 328 & 10 & $L_{1}^{4}L_{3}^{6}L_{5}^{4}L_{6}p_1^{6}$ & 54,526 & 9 & 70\\
6 & $L_{1}L_{3}L_{4}L_{5}L_{6}p_1^{2}$ & 54,448 & 164 & 1000 & 10 & $L_{1}L_{3}L_{4}L_{5}L_{6}^{2}p_1^{6}$ & 58,582 & 14 & 905\\
6 & $L_{1}L_{3}L_{4}L_{5}L_{6}p_1^{2}$ & 54,454 & 150 & 474 & 10 & $L_{1}^{3}L_{2}L_{3}^{4}L_{5}^{3}L_{6}p_1^{5}$ & 54,536 & 38 & 186\\
6 & $L_{1}^{3}L_{3}^{3}L_{7}\mathsf{D}_{4}^{2}p_1$ & 56,486 & 396 & 740 & 10 & $L_{1}^{3}L_{3}^{4}L_{4}L_{5}^{3}L_{6}p_1^{5}$ & 54,524 & 38 & 905\\
6 & $L_{1}^{2}L_{2}L_{3}\mathsf{D}_{4}^{2}p_1$ & 52,400 & 933 & 207 & 11 & $L_{4}L_{5}^{2}p_1^{9}$ & 50,468 & 3 & 30\\
6 & $L_{1}^{2}L_{2}L_{3}L_{5}\mathsf{D}_{4}p_1$ & 54,440 & 396 & 492 & 11 & $L_{1}L_{3}L_{4}L_{5}L_{6}p_1^{8}$ & 54,524 & 8 & 40\\
6 & $L_{1}L_{2}L_{3}^{2}L_{5}\mathsf{D}_{4}p_1$ & 54,446 & 352 & 450 & 11 & $L_{3}^{2}L_{4}L_{6}^{2}p_1^{8}$ & 58,580 & 6 & 110\\
6 & $L_{1}^{2}L_{2}L_{3}L_{5}L_{6}p_1$ & 54,460 & 228 & 2420 & 11 & $L_{1}^{2}L_{3}L_{4}L_{5}L_{6}^{2}p_1^{7}$ & 58,590 & 10 & 7\\
6 & $L_{1}^{2}L_{3}L_{4}L_{5}L_{6}p_1$ & 54,448 & 228 & 279 & 11 & $L_{1}^{4}L_{2}L_{3}^{6}L_{5}^{4}L_{6}p_1^{6}$ & 54,560 & 9 & 186\\
6 & $L_{1}L_{2}L_{3}^{2}L_{5}L_{6}p_1$ & 54,466 & 128 & 1490 & 11 & $L_{1}^{4}L_{3}^{6}L_{4}L_{5}^{4}L_{6}p_1^{6}$ & 54,548 & 9 & 186\\
6 & $L_{1}L_{2}L_{3}L_{5}^{2}L_{6}p_1$ & 54,476 & 178 & 628 & 12 & $L_{1}L_{4}L_{5}^{2}L_{6}p_1^{9}$ & 54,548 & 4 & 49\\
6 & $L_{1}L_{3}^{2}L_{4}L_{5}L_{6}p_1$ & 54,454 & 128 & 328 & 12 & $L_{1}L_{3}^{2}L_{4}L_{5}L_{6}^{2}p_1^{8}$ & 58,604 & 7 & 20\\
6 & $L_{1}L_{3}L_{4}L_{5}^{2}L_{6}p_1$ & 54,464 & 178 & 474 & 13 & $L_{1}^{2}L_{3}^{2}L_{4}L_{5}^{2}L_{6}^{2}p_1^{9}$ & 58,628 & 4 & 55\\
7 & $L_{1}L_{3}^{2}L_{5}L_{6}p_1^{4}$ & 54,454 & 64 & 92 & 13 & $L_{1}^{2}L_{3}L_{4}L_{5}^{2}L_{6}^{2}p_1^{9}$ & 58,628 & 3 & 27\\
7 & $L_{2}L_{3}^{2}L_{6}p_1^{4}$ & 54,464 & 34 & 320 & 15 & $L_{3}^{3}L_{4}L_{6}^{3}p_1^{12}$ & 62,708 & 1 & 4\\
7 & $L_{3}^{2}L_{4}L_{6}p_1^{4}$ & 54,452 & 27 & 72 &  & & & & \\
\cline{1-10}
\end{tabular}
\end{center}
\end{scriptsize}
\end{table}

\end{document}